\documentstyle[12pt,leqno]{article}
\begin{document}
\title{A Survey on Nambu-Poisson Brackets}
\author{{\normalsize by}\\Izu Vaisman}
\date{}
\maketitle
{\def\thefootnote{*}\footnotetext[1]%
{{\it 1991 Mathematics Subject Classification}
58 F 05. \newline\indent{\it Key words and phrases}:
Nambu-Poisson brackets, Nambu-Poisson tensors, Nambu-Poisson-Lie groups.}}
\begin{center} \begin{minipage}{12cm}
A{\footnotesize BSTRACT.
The paper provides a survey of known results
on geometric aspects related to Nambu-Poisson brackets.}
\end{minipage} \end{center}
\section{Introduction}
In 1973, Nambu \cite{Nb} studied a dynamical system which was defined as a
Hamiltonian system with respect to
a ternary Poisson bracket.
A few other papers on this bracket have followed at the time
\cite{BF}, \cite{MS}. A few years ago, Takhtajan \cite{Tk1} reconsidered the
subject, proposed a general, algebraic definition of a {\it
Nambu-Poisson bracket of order} $n$, and gave the basic characteristic
properties of this operation. The Nambu-Poisson bracket
is an intriguing operation, in spite of its rather restrictive
character, which follows from the fact
conjectured in \cite{Tk1}, and proven by several authors
\cite{Wz} (cited by \cite{DFS}, and much older than \cite{Tk1}),
\cite{Gt}, \cite{AG},
\cite{Nks}, \cite{Le}, \cite{MVV} namely, that,
locally and with respect to well chosen coordinates,
any Nambu-Poisson bracket is just a Jacobian
determinant  as in \cite{Nb}. In particular,
the deformation quantization of the Nambu-Poisson bracket
leads to interesting mathematical developments
\cite{DFS}, \cite{DF}. On the other hand, the bracket inspired some
generalizations of Lie-algebraic constructions (anticipated in \cite{Fp})
\cite{ChT}, \cite{Tk2},
\cite{TkD}, \cite{Gt}, \cite{AP}, \cite{MV}.

The aim of this paper is to give a survey of the
subject from the point of view of geometry. In the next section, we
review the basics, and present the geometric structure of Nambu-Poisson
manifolds. Another section will be devoted to Nambu-Poisson-Lie groups.
Finally, while
we do not intend to review quantization theories, we formulate some related
questions in the last section.

The paper does not contain new results.
Everything in the paper is in the $C^{\infty}$ category. Information on
the usual Poisson manifolds may be found in \cite{V1}, for instance.
More general Nambu-Jacobi brackets were also studied 
\cite{MarJ}, \cite{{MVV},
{GM}} but, we will not discuss this subject here.

{\it Acknowledgements}.
The author thanks N. Nakanishi, J. C. Marrero,
and J.-P. Dufour for useful remarks.
The final version was written during a stay at the Erwin
Schr\"odinger Institute for Mathematical Physics in Vienna, Austria. I want to
express here my gratitude to the ESI for its support, and to Prof. Peter Michor for his
invitation there and for our discussions on the subject.
\section{Nambu-Poisson Brackets}
Let $M^{m}$ be an $m$-dimensional differentiable manifold, and ${\cal
F}(M)$ its algebra of real valued $C^{\infty}$-functions. A {\it
Nambu-Poisson bracket or structure} of order $n$, $3\leq n\leq m$
(this condition is always imposed in the paper) is an internal
$n$-ary operation on ${\cal F}(M)$, denoted by $\{\;\}$, which satisfies the
following axioms:\\
$(i)$ $\{\;\}$ is {\bf R}-multilinear and totally skew-symmetric;\\
$(ii)$ \hspace{1cm}$\{f_{1},\ldots,f_{n-1},gh\}=\{f_{1},\ldots,f_{n-1},g\}h
+g\{f_{1},\ldots,f_{n-1},h\}$\\
(the {\it Leibniz rule});\\
$(iii). \hspace{2cm}\{f_{1},\ldots,f_{n-1},\{g_{1},\ldots,g_{n}\}\}=$
$$\sum_{k=1}^{n}\{g_{1},\ldots,g_{k-1},
\{f_{1},\ldots,f_{n-1},g_{k}\},g_{k+1},\ldots,g_{n}\}$$
(the {\it fundamental identity}).
A manifold endowed with a Nambu-Poisson bracket is a {\it Nambu-Poisson
manifold}. Remember that if we use the same definition for $n=2$, we get a
{\it Poisson bracket}.

By $(ii)$, $\{\;\}$ acts on each factor as a vector field, whence it must be
of the form
$$\{f_{1},\ldots,f_{n}\}=P(df_{1},\ldots,df_{n}),\leqno{(2.1)}$$
where $P$ is a field of $n$-vectors on $M$. If such a field defines a
Nambu-Poisson bracket, it is called a {\it Nambu-Poisson tensor (field)}.
$P$ defines a bundle mapping
$$\sharp_{P}:\underbrace
{T^{*}M\times\ldots\times T^{*}M}_{n-1\:times}
\longrightarrow TM \leqno{(2.2)}$$
given by
$$<\beta,\sharp_{P}(\alpha_{1},\ldots,\alpha_{n-1})>=
P(\alpha_{1},\ldots,\alpha_{n-1},\beta) \leqno{(2.3)}$$
where all the arguments are covectors.

In what follows, we denote an $n$-sequence of functions or forms, say
$f_{1},\ldots, f_{n}$, by $f_{(n)}$, and, if an index $k$ is missing, by
$f_{(n,\hat k)}$.

The next basic notion is that of the $P$-{\it Hamiltonian vector
field} of $(n-1)$ functions defined by
$$X_{f_{(n-1)}}=\sharp_{P}(df_{(n-1)}).\leqno{(2.4)}$$
Then, the fundamental identity $(iii)$ means that {\it the Hamiltonian vector
fields are derivations of the Nambu-Poisson bracket}.

Another
interpretation of $(iii)$ is
$$(L_{X_{f_{(n-1)}}}P)(dg_{1},\ldots,dg_{n})=0,\leqno{(2.5)}$$
where $L$ is the Lie derivative, i.e., {\it the Hamiltonian vector fields
are infinitesimal automorphisms of the Nambu-Poisson tensor}.

The
fundamental identity also implies
$$X_{f_{(n-1)}}X_{g_{(n-1)}}h=
\sum_{k=1}^{n-1}\{g_{1},\ldots,g_{k-1},
X_{f_{(n-1)}}g_{k},g_{k+1},\ldots,g_{n-1},h\}
+X_{g_{(n-1)}}X_{f_{(n-1)}}h,$$
whence $$[X_{f_{(n-1)}},X_{g_{(n-1)}}]=
\sum_{k=1}^{n-1}X_{(g_{1},\ldots,g_{k-1},
X_{f_{(n-1)}}g_{k},g_{k+1},\ldots,g_{n-1})}.\leqno{(2.6)}$$
Therefore, {\it the set ${\cal H}(P)$ of all the real,
finite, linear combinations
of Hamiltonian vector fields is a Lie algebra}. (Notice that for $n\geq3$
such a combination may not be a Hamiltonian vector field itself!)

The Nambu-Poisson tensor fields were characterized as follows by Takhtajan
\cite{Tk1} \proclaim 2.1 Theorem. The $n$-vector field
$P$ is a Nambu-Poisson tensor of order $n$
($n\geq3$) iff the natural components of $P$ with respect to any local
coordinate system $x^{a}$ of $M$ satisfy the equalities:
$$\sum_{k=1}^{n}[P^{b_{1}\ldots b_{k-1}ub_{k+1}\ldots b_{n}}
P^{va_{2}\ldots a_{n-1}b_{k}}+
P^{b_{1}\ldots b_{k-1}vb_{k+1}\ldots b_{n}}
P^{ua_{2}\ldots a_{n-1}b_{k}}]=0,\leqno{(2.7)}$$
$$\sum_{u=1}^{m}[P^{a_{1}\ldots a_{n-1}u}\partial_{u}P^{b_{1}\ldots
b_{n}}
-\sum_{k=1}^{n}P^{b_{1}\ldots b_{k-1}ub_{k+1}\ldots b_{n}}
\partial_{u}P^{a_{1}\ldots a_{n-1}b_{k}}]=0.\leqno{(2.8)}$$
Furthermore, $P$ is a Nambu-Poisson tensor field iff $P/_{U}$ 
is a Nambu-Poisson tensor field, for
$U=\{x\in M\:/\:P_{x}\neq0\}$.
\par
\noindent{\bf Proof.} Fix a point $p\in M$, and local coordinates $x^{i}$
around $p$ such that $x^{i}(p)=0$. Then,
with the Einstein summation convention, denote
$$P=\frac{1}{n!}P^{i_{1}\ldots i_{n}}
\frac{\partial}{\partial x^{i_{1}}}\wedge\ldots\wedge
\frac{\partial}{\partial x^{i_{n}}},\leqno{(2.9)}$$
$$\partial_{u}P=\frac{1}{n!}\frac{\partial P^{i_{1}\ldots i_{n}}}
{\partial x^{u}}\frac{\partial}{\partial x^{i_{1}}}\wedge\ldots\wedge
\frac{\partial}{\partial x^{i_{n}}},\leqno{(2.10)}$$
$$\partial P=\partial_{u}P\otimes dx^{u}.\leqno{(2.11)}$$
If the fundamental identity is expressed by means of 
(2.1), the terms which contain the second derivatives of the functions $g$
cancel, and the identity becomes
$$\sum_{u=1}^{m}P(df_{(n-1)},dx^{u})(\partial_{u}P)(dg_{(n)})
\leqno{(2.12)}$$
$$=\sum_{u=1}^{m}\sum_{k=1}^{n}[P(dg_{1},\ldots,dg_{k-1},dx^{u},
dg_{k+1},\ldots,dg_{n})(\partial_{u}P)(df_{(n-1)},dg_{k})$$
$$+\sum_{h=1}^{n-1}P(dg_{1},\ldots,dg_{k-1},dx^{s},
dg_{k+1},\ldots,dg_{n})P(df_{1},\ldots,df_{h-1},$$
$$\frac{\partial^{2}f_{h}}
{\partial x^{s}\partial x^{t}}dx^{t},
df_{h+1},\ldots,df_{n-1},dg_{k})].$$

Now, (2.12) is always true at $p$ if it is true in the following two
cases:\\
a) $f_{i}=x^{a_{i}},\;g_{j}=x^{b_{j}}$,
b) same as in a) with the exception of $f_{1}=x^{u}x^{v}$. Case a) yields
(2.8), and case b) yields (2.7).
Finally, the
restriction $x^{i}(p)=0$ may be removed by a translation of the
coordinates.

The last assertion of the theorem is an obvious consequence of (2.7), (2.8).
Q.e.d.

Equality (2.7) is algebraic, and
it is called the {\it quadratic identity}. This
condition does not appear for the usual Poisson structures ($n=2$).
Equality (2.8) is called the
{\it differential identity}, and it does not have a
tensorial character. However, it is clear that if (2.7), (2.8) hold for one
coordinate system at $p\in M$  the fundamental identity holds hence, (2.7),
(2.8) will hold in any coordinate system.

The quadratic identity is rather intriguing. For this reason, we give several
equivalent expressions below. First, (2.7) is equivalent with
$$\sum_{k=1}^{n}[\{\varphi,f_{1},\ldots,f_{n-2},g_{k}\}
\{\psi,g_{1},\ldots,\hat g_{k},\ldots,g_{n}\}\leqno{(2.13)}$$
$$+\{\psi,f_{1},\ldots,f_{n-2},g_{k}\}
\{\varphi,g_{1},\ldots,\hat g_{k},\ldots,g_{n}\}]$$
for arbitrary functions. Indeed, using (2.1) we see that (2.7) implies (2.13), and
on
the other hand
(2.13) reduces to (2.7) in the case of the
coordinate functions.

Then, the expression (2.12) of the fundamental identity
is the same as
$$<\sharp_{P}(df_{(n-1)}),\partial P(dg_{(n)})>\leqno{(2.14)}$$
$$=\sum_{k=1}^{n}(-1)^{n-k}
[<\sharp_{P}(dg_{(n,\hat k)}),\partial P(df_{(n-1)},dg_{k})>$$
$$+\sum_{h=1}^{n-1}(-1)^{h+k}(Hess\,f_{h})(\sharp_{P}(dg_{(n,\hat k)}),
\sharp_{P}(df_{(n-1),\hat h)},dg_{k}))],$$
where all $f\in{\cal F}(M)$, and
$$Hess\,f:=\frac{\partial^{2}f}
{\partial x^{s}\partial x^{t}}dx^{s}\otimes dx^{t}$$
is the non invariant Hessian of $f$.

Moreover, if $\nabla$ is an arbitrary torsionless connection on $M$, (2.14)
is equivalent with the same relation where the partial derivatives in
$\partial P$ and in the Hessians are replaced by $\nabla$-covariant
derivatives. This yields a tensorial expression of the
fundamental identity.

Finally, (2.14) yields another invariant expression
of the quadratic identity if we proceed as follows. Notice that the quadratic identity
holds iff (2.14) holds
for functions which have a vanishing second derivatives
at the point $p$, except for
$f_{1}$, for which we ask the vanishing of the first derivatives, while
$Hess\,f_{1}=T$ is an arbitrary $2$-covariant symmetric tensor.
Accordingly, the quadratic identity is equivalent to
$$\sum_{k=1}^{n}(-1)^{k+1}T(\sharp_{P}(\lambda_{(n,\hat k)}),
\sharp_{P}(\mu_{(n-1,\hat 1)},\lambda_{k}))=0\leqno{(2.15)}$$
for any $2$-covariant, symmetric tensor $T$, and any covectors
$\lambda,\mu$.

The geometric meaning of the quadratic identity will
be shown in the forthcomming Theorem 2.4.

A mapping $\varphi:(M_{1},P_{1})
\rightarrow(M_{2},P_{2})$
between two Nambu-Poisson manifolds of the same order $n$
is a {\it Nambu-Poisson morphism} if the tensor
fields $P_{1}$ and $P_{2}$
are
$\varphi$-related or, equivalently,
$\forall g_{(n)}\in{\cal F}(M_{2})$, one has
$$\{g_{1}\circ\varphi,\ldots,g_{n}\circ\varphi\}_{1}=
\{g_{1},\ldots,g_{n}\}_{2}.$$
Moreover, if $\varphi$ is a diffeomorphism, the two manifolds are said to
be {\it equivalent Nambu-Poisson manifolds}. The notion of a 
Nambu-Poisson morphism also allows us to give the following definition:
a submanifold $N$ of the Nambu-Poisson manifold $(M,P)$
is a Nambu-Poisson submanifold if $N$ has a (necessarily unique)
Nambu-Poisson tensor field $Q$ of the same order as $P$
such that the inclusion of $(N,Q)$ in $(M,P)$ is a Nambu-Poisson
morphism. As in the Poisson case $n=2$, $Q$ exists iff, along $N$,
$P$ vanishes whenever evaluated on $n$ $1$-forms one of which, at least, belongs to 
the annihilator space $Ann(TN)$,
and then $im\sharp_{P}$ is a tangent distribution of $N$
e.g., \cite{V1}.

By Theorem 2.1 $P$ is a Nambu-Poisson tensor on the manifold $M$
iff it is such on its nonvanishing subset. The
following theorem \cite{{Gt},{AG},{Nks},{Le},{MVV}} establishes
the local canonical
structure of the Nambu-Poisson brackets around nonvanishing points, up to equivalence.
\proclaim 2.2 Theorem.
$P$ is a Nambu-Poisson tensor field of order $n$ iff $\forall p\in M$
where $P_{p}\neq0$ there are local coordinates
$(x^{k},y^{\alpha})$ ($k=1,\ldots,n$, $\alpha=1,\ldots,m-n$) around $p$
such that $$P=\frac{\partial}{\partial x^{1}}\wedge\ldots\wedge
\frac{\partial}{\partial x^{n}}\leqno{(2.16)}$$
on the corresponding coordinate neighborhood.\par
\noindent{\bf Proof.} If (2.16) holds, we have $P^{1...n}=1$, and the
components of $P$ which have other indices than a permutation of $(1,...,n)$
vanish. It is easy to see that (2.7), (2.8) hold in this case.

The following proof of the converse result
belongs to Nakanishi \cite{Nks},
and is modeled on Weinstein's proof of the local structure theorem of
Poisson manifolds (e.g., \cite{W}, \cite{V1}).
Around $p$, take functions $x_{(n-1)}$ such that $X_{x_{(n-1)}}\neq0$, then
change to local coordinates $z_{(m)}$ where $X_{x_{(n-1)}}=\partial/
\partial z_{1}$, and put $x_{n}=z_{1}$.
Since
$$\{x_{1},\ldots,x_{n}\}=1, \leqno{(2.17)}$$
$x_{(n)}$ are functionally independent, and the vector fields $Y_{k}:=
(-1)^{n-k}X_{x_{(n,\hat k)}}$,
which satisfy $Y_{k}(x_{h})=\delta_{kh}$, are
linearly independent. Moreover, (2.6) shows that $Y_{k}$ commute, and there
exist local coordinates $(s_{k},y_{\alpha})$
($k=1,\ldots,n$, $\alpha=1,\ldots,m-n$)
such that $Y_{k}=\partial/\partial s_{k}$ for all
$k$. Furthermore, by
looking at the corresponding Jacobian, we see that $(x_{k},
y_{\alpha})$ also are local coordinates around $p$, and such that
$Y_{k}=\partial/\partial x_{k}$, and all $\{x_{k_{1}},\ldots,x_{k_{n-1}},
y_{\alpha}\}
=0$.

The following trick is to evaluate in two ways the bracket
$$\frac{1}{2}(-1)^{k-1}\{x_{1}^{2},x_{2},\ldots,x_{n-1},
\{x_{2},\ldots,x_{k},x_{n},y_{\alpha_{1}},\ldots,y_{\alpha_{h}}\}\}$$
where $k+h=n$. If we use first the fundamental identity and then the
Leibniz rule we get $\{x_{1},\ldots,x_{k},
y_{\alpha_{1}},\ldots,y_{\alpha_{h}}\}$. If we use first the Leibniz rule
and then the fundamental identity, we get $0$.
(Use (2.17) in both computations.) Similarly, we get the
general result
$$P^{i_{1}\ldots i_{k}\alpha_{1}\ldots\alpha_{h}}
=\{x_{i_{1}},\ldots,x_{i_{k}},
y_{\alpha_{1}},\ldots,y_{\alpha_{h}}\}=0. \leqno{(2.18)}$$

Finally, we must compute the components of $P$ with Greek indices only. Of
course, they vanish if $m<2n$. If $m\geq2n\geq6$, these components are again
given by using (2.17), (2.18) and a two-way
computation of a Nambu bracket namely,
$$0=\{x_{1}y_{\alpha_{1}},x_{2},\ldots,x_{n-1},\{x_{n},y_{\alpha_{2}},
\ldots,y_{\alpha_{n}}\}\}\leqno{(2.19)}$$
$$=\{y_{\alpha_{1}},y_{\alpha_{2}},
\ldots,y_{\alpha_{n}}\}=
P^{\alpha_{1}\ldots\alpha_{n}}.$$
The results (2.17), (2.18), (2.19), with the notational change of writing the
indices of the coordinates up as usual,
imply (2.16). Q.e.d.
\proclaim 2.3 Remark. On the canonical coordinate neighborhood where
(2.16) holds we have $${\cal
D}:=span(im\,\sharp_{P})=span\{\partial/\partial x^{k}\}.$$
Hence, globally ${\cal D}$
is a foliation with singularities whose leaves are either points,
called {\it singular points} of $P$,
or $n$-dimensional submanifolds with a Nambu-Poisson
bracket induced by $P$. (In other words, the
computation of the latter is along the leaves of ${\cal D}$).\par
This remark extends well
known results of Poisson geometry (e.g., \cite{V1}), and
it was proven in \cite{Gt} and
\cite{Le}. In \cite{Le} the proof is
by applying the Stefan-Sussmann-Frobenius
theorem to  ${\cal D}$, which is possible because ${\cal D}$ is also equal to
$span\,{\cal H}(P)$. We call ${\cal D}$ the
{\it canonical foliation} of the Nambu-Poisson structure $P$. The canonical
foliation is regular i.e., all the leaves are $n$-dimensional, iff
$P$ never vanishes, and then we will say that $P$ is a {\it regular
Nambu-Poisson structure}. \vspace{2mm}

The structure theorem 2.2 allows us to prove one more important result.
First, we will say that an $n$-vector field is {\it decomposable} (or {\it
simple}) if, $\forall p\in M$, there are $V_{1},...,V_{n}\in T_{p}M$ such
that $P_{p}=V_{1}\wedge...\wedge V_{n}$. (This does not mean that such a
decomposition holds for global vector fields on $M$.) Then, we
have
\proclaim 2.4 Theorem. The quadratic identity (2.7) is equivalent with
the fact that  the $n$-vector field $P$ is decomposable. \par
\noindent{\bf Proof.} This is a pointwise, algebraic result, and an
algebraic proof can be found in \cite{AG}. On the other hand
in \cite{DFS} the authors quote
\cite{Wz} for a proof of the result which is much older than the notion of
a Nambu-Poisson bracket.

Here, we will use Theorem 2.2. Clearly, it suffices to prove the result on
${\bf R}^{m}$.
If $P$ is decomposable, we use a vector basis which has
$V_{1},...,V_{n}$ as its first vectors, and a straightforward inspection of
(2.7) shows that this condition holds.

Conversely, if $P$ is given at a point, and it satisfies the quadratic
identity, we may extend it to a tensor
field with constant components on ${\bf
R}^{m}$. The latter then obviously also satisfies the differential identity
(2.8), and is a Nambu-Poisson tensor field on ${\bf R}^{m}$. Thus, $P$ is
decomposable by Theorem 2.2. Q.e.d.\vspace{2mm}

In connection with Theorem 2.4, let us remember that decomposable $n$-vectors are also
characterized by the {\it Pl\"ucker relations} (e.g., \cite{Shaf}, p.42)
$$(-1)^{n}P^{uva_{2}...a_{n-1}}P^{b_{1}...b_{n}}=
\sum_{k=1}^{n}P^{ua_{2}...a_{n-1}b_{k}}P^{b_{1}...b_{k-1}vb_{k+1}...b_{n}}.$$
By a symmetrization, these relations yield (2.7), and Theorem 2.4 tells us
that (2.7) are equivalent to the Pl\"ucker relations.

Another immediate consequence of Theorem 2.2 is \cite{Le}
\proclaim 2.5 Corollary. A
Nambu-Poisson tensor field $P$ of an even order $n=2s$
satisfies the condition $[P,P]=0$, where the operation is the
Schouten-Nijenhuis bracket.
\par
This corollary suggests the study of {\it generalized Poisson structures}
\cite{AP}, \cite{Le} defined by a $(2s)$-vector field $P$ such that
$$[P,P]=0.\leqno{(2.20)}$$

The canonical expression (2.16) provides
the basic example of a Nambu-Poisson bracket,
which was considered in Nambu's original paper \cite{Nb} for $n=3$.
Namely, (2.16) means that we have
$$\{f_{1},\ldots,f_{n}\}=\frac{\partial
(f_{1},\ldots,f_{n})}{\partial
(x^{1},\ldots,x^{n})}.\leqno{(2.21)}$$

This example may be extended to a description of all the regular
Nambu-Poisson structures
\cite{Gt}, \cite{Le}.
\proclaim 2.6 Theorem. A regular Nambu-Poisson structure of
order $n$ on a
differentiable manifold $M^{m}$is the same thing as a regular $n$-dimensional
foliation $S$ of $M$, and a bracket operation defined by the formula
$$d_{S}f_{1}\wedge\ldots\wedge d_{S}f_{n}=\{f_{1},\ldots,f_{n}\}\omega,
\leqno{(2.22)}$$
where $\omega$ is an $S$-leafwise volume form,
and $d_{S}$ is differentiation along the leaves of $S$. \par
\noindent{\bf Proof.}
First, let $M^{m}$ be a differentiable manifold endowed with a regular
$n$-dimensional foliation $S$, and an $S$-leafwise volume form $\omega$.
(E.g., see \cite{Mol} for foliation theory.)
Then, the bracket defined by (2.22) is
a regular Nambu-Poisson bracket.
Indeed, if $x^{(n)}$ are local coordinates along the
leaves
of $S$, and if $$\omega=\varphi d_{S}x^{1}\wedge\ldots\wedge d_{S}x^{n},$$
we get the local expression
$$\{f_{1},\ldots,f_{n}\}=\frac{1}{\varphi}\frac{\partial
(f_{1},\ldots,f_{n})}{\partial
(x^{1},\ldots,x^{n})}.\leqno{(2.23)}$$

Then, the change of the local coordinates
$$\tilde x^{1}=\int\varphi dx^{1},\;\tilde x^{2}=x^{2},\ldots,x^{n}=x^{n}$$
leads to (2.21) in the new coordinates $\tilde x^{(n)}$.

In particular, notice from the proof above
that any formula of the type (2.23) defines a
regular Nambu-Poisson bracket.

Now, conversely, if $P$ is a regular Nambu-Poisson structure on $M$, we
take $S$ to be the canonical foliation of $P$, and choose the
leafwise volume form $\omega$ such that $i(P)\omega=1$. Then, we see that
(2.22) holds by applying to it the operator $i(P)$.
Clearly, the chosen volume form is the only possible one. Q.e.d.
\vspace{2mm}

Following is a number of other interesting facts relevant to Nambu-Poisson structures.
\proclaim 2.7 Remarks. i) \cite{GM}. A decomposable $n$-vector field $P$
is a Nambu-Poisson tensor iff the distribution ${\cal D}=span(im\sharp_{P})$
is involutive on the set of the non singular points of $P$.
ii) \cite{Tk1}. If
we have a Nambu-Poisson bracket of order $n>2$, and keep $p$ of its arguments
fixed, we get a Nambu-Poisson bracket of order $n-p$ (a Poisson bracket if
$n-p=2$), and, conversely 
\cite{GM}, if the result of an arbitrary fixed choice of $p$ arguments
($p=1,...,n-2$) always yields a Nambu-Poisson tensor, $P$ is a 
Nambu-Poisson tensor.
iii). If $(M_{a},P_{a})$  are Nambu-Poisson manifolds of
order $n_{a}\geq3$ ($a=1,2$),
then $(M_{1}\times M_{2},P_{1}\wedge P_{2})$ is a
Nambu-Poisson manifold of order $n_{1}+n_{2}$.
iv). If $P$ is a Nambu-Poisson tensor on a manifold $M$, so is $fP$ for
any function $f\in C^{\infty}(M)$. In particular, this implies that (2.7)
is equivalent to
$$P^{i_{1}...i_{n-1}k}P^{j_{1}...j_{n}}=\sum_{h=1}^{n}P^{j_{1}...j_{h-1}
kj_{h+1}...j_{n}}P^{i_{1}...i{n-1}j_{h}}. \leqno{(2.7')}$$
\par
Concerning the first remark, we already know that ${\cal D}$ is
involutive whenever $P$ is Nambu-Poisson (Remark 2.3). On the other hand,
since for $P=V_{1}\wedge...\wedge V_{n}$, ${\cal D}=span\{V_{1},...,V_{n}\}$,
if ${\cal D}$ is involutive, we have $P=(\partial/\partial x_{1})
\wedge...\wedge(\partial/\partial x_{n})$ in some well chosen local coordinates 
on a neighborhood of $x\in M$ 
where $P_{x}\neq 0$ (Frobenious Theorem). Then, the corresponding
bracket takes the form (2.23), and it is a Nambu-Poisson bracket.

The direct part of the second remark follows by
checking the axioms. For the converse, it suffices to take $p=1$, and 
check by a computation
that if (2.7), (2.8) hold for $i(df)P$, $\forall f\in C^{\infty}(M)$, 
they also hold for $P$ itself.

The third remark is an immediate consequence of (2.16).

The last remark follows by putting $P$ under the form (2.16), and
using the proof of Theorem 2.6. Then, (2.7$'$) is the coordinate expression
of the fact that $fP$ satisfies the fundamental identity
$\forall f\in C^{\infty}(M)$. (It is obvious
that (2.7$'$) implies (2.7).) For arbitrary functions (2.7$'$) yields
$$\{f_{1},...,f_{n-1},f\}\{g_{1},...,g_{n}\}\leqno{(2.7'')}$$
$$=\sum_{h=1}^{n}\{g_{1},...,g_{h-1},f,g_{h+1},...,g_{n}\}\{f_{1},...,f_{n-1},
g_{h}\}.$$

The structure theorem 2.2 was used by Dufour and Zung
\cite{DZ}, and by Nakanishi \cite{Nak} in order to
characterize Nambu-Poisson manifolds by means of differential forms, which
are better suited for calculus than the multivectors. Namely, if $\omega$
is a volume form on the manifold $M^{m}$, for every $n$-vector $P$ there
exists a corresponding $(m-n)$-form $\varpi:=i(P)\omega$, and the result
proven in \cite{DZ} is that $P$ is a Nambu-Poisson tensor iff
$$(i(A)\varpi)\wedge\varpi=0,\hspace{3mm} (i(A)\varpi)\wedge d\varpi=0,
\leqno{(2.24)}$$
for any $(m-n-1)$-vector $A$.
In \cite{DZ}, a differential form $\varpi$ which satisfies
(2.24) is called a {\it Nambu co-form}.
In \cite{Nak} it is shown that $\varpi$ is a Nambu co-form iff it is
decomposable and integrable i.e., $d\varpi=\theta\wedge\varpi$ for some
$1$-form $\theta$.
\vspace{2mm}

On ${\bf R}^{m}$, any constant, decomposable $n$-vector field
$k^{i_{1}\ldots i_{n}}$ is a Nambu-Poisson tensor, since it satisfies both
the quadratic and the differential identities. If we use Remark 2.7 ii) for
this Nambu-Poisson tensor $k$, and keep as a
fixed function $(1/2)\sum_{j=1}^{m}
(x^{j})^{2}$, we get a new Nambu-Poisson tensor, of order $n-1$,
with the natural components
$$P^{i_{1}\ldots i_{n-1}}=\sum_{j=1}^{m}k^{i_{1}\ldots i_{n-1}j}x^{j}.
\leqno{(2.25)}$$

A Nambu-Poisson structure defined on ${\bf R}^{m}$ by a tensor whose
natural components are linear functions of $x^{j}$ is called a
{\it linear Nambu-Poisson structure}, and (2.25) gives the basic example
\cite{ChT}. Linear Nambu-Poisson
structures are a generalization of the Lie-Poisson structures of Lie
coalgebras
(e.g., \cite{V1}).
Accordingly, a definition and study of {\it $n$-Lie algebras}
is suggested \cite{{Fp}, {Tk1}, {Tk2}, {TkD}, {Gt}, {MV}, {MVV}}.
More precisely, a $n$-Lie algebra (called {\it Fillipov algebra} in \cite{GM})
is a vector space endowed with an
internal, $n$-ary, skew symmetric bracket which satisfies the fundamental
identity of a Nambu-Poisson bracket.
(Different notions of $n$-Lie algebras were studied in \cite{HW} and
\cite{MV}.)
By looking at brackets of linear functions, it
easily follows  that a linear Nambu-Poisson structure of order $n$ on
${\bf R}^{m}$ induces a $n$-Lie algebra structure on the dual of
${\bf R}^{m}$ \cite{Tk1}. The converse may not be true since the structure
constants of a general $n$-Lie algebra may form a non decomposable $n$-vector.

For instance, if $m=n+1$ we may take $k$ in (2.25) to be the canonical volume
tensor of ${\bf R}^{n+1}$, and we get the linear Nambu-Poisson structure of
order $n$ discussed in \cite{ChT}. The corresponding $n$-Lie algebra is the
vector space ${\bf R}^{n+1}$ endowed with the operation of a
{\it vector product}
of $n$ vectors (the determinant which has the coordinates of the vectors,
and
the canonical, positive, orthonormal basis as its columns
\cite{Bg}). Another definition of this operation, denoted by $\times$, is
$$v_{1}\times\ldots\times v_{n}=*(v_{1}\wedge\ldots\wedge v_{n}),
\leqno{(2.26)}$$
where $*$ is the Hodge star operator of the canonical Euclidean metric of
${\bf R}^{n+1}$.
It is also easy to see that the
canonical foliation of the linear Nambu-Poisson structure of
${\bf R}^{n+1}$ defined above has the origin as a $0$-dimensional leaf, and
the spheres with center at the origin as $n$-dimensional leaves. (For
$n=2$,
this is the dual of the Lie algebra $o(3)$ with its well known Lie-Poisson
structure.)\vspace{2mm}

Of course, we may replace ${\bf R}^{m}$ by any vector space , with linear
coordinates, in the definition of a linear Nambu-Poisson structure. Then,
as in the case of Poisson structures \cite{W},
we notice that, if $(M,P)$ is a
Nambu-Poisson manifold, and if $p\in M$ is a {\it singular point} of $P$
(i.e., $P(p)=0$), the linear part of the Taylor development of $P$ at $p$
defines a linear Nambu-Poisson structure on $T_{p}M$, and a corresponding
$n$-Lie algebra structure on $T_{p}^{*}M$, which are independent of the
choice of the local coordinates at $p$. This
linear Nambu-Poisson structure of $T_{p}M$ should be regarded as the
{\it linear approximation} of $P$ at $p$, and $P$ is {\it linearizable} at
$p$ if $P$ is equivalent with its linear approximation on some
neighbourhood of $p$.

The linear Nambu-Poisson tensors are completely determined by Dufour and
Zung in \cite{DZ} (see also \cite{MVV} and \cite{GM}), 
and the result is \proclaim 2.8 Theorem.
For any linear Nambu-Poisson structure $P$ of order $n$ on the linear space
$V^{m}$ there exists a basis of $V$ such that the tensor
$P$ is of one of the following types.\\
Type I:
$$P=\sum_{j=1}^{r+1}\pm x_{j}\frac{\partial}{\partial x_{1}}\wedge...
\wedge\frac{\partial}{\partial x_{j-1}}\wedge
\frac{\partial}{\partial x_{j+1}}\wedge...\wedge
\frac{\partial}{\partial x_{n+1}}$$
$$+\sum_{j=1}^{s}\pm x_{n+j+1}\frac{\partial}{\partial
x_{1}}\wedge...\wedge
\frac{\partial}{\partial x_{r+j}}\wedge
\frac{\partial}{\partial x_{r+j+2}}\wedge...\wedge
\frac{\partial}{\partial x_{n+1}},$$
with $-1\leq r\leq n,\;0\leq s\leq min(m-n-1,n-r)$;\\
Type II:
$$P=\frac{\partial}{\partial x_{1}}\wedge...\wedge
\frac{\partial}{\partial x_{n-1}}\wedge(\sum_{i,j=n}^{m}a^{i}_{j}x_{i}
\frac{\partial}{\partial x_{j}}).$$
\par
In the proof of Theorem 2.8 an essential role is played by the following results of
linear geometry (Lemma 3.2 and Theorem 3.1 of \cite{DZ}, Lemma 1 of \cite{GM})
\proclaim 2.9 Lemma. i). Let $P_{1},P_{2}$ be decomposable $n$-vectors of a linear space
$V$ such that $P_{1}+P_{2}$ is also decomposable and let $D_{1},D_{2}$ be the
subspaces
spaned by the factors of $P_{1},P_{2}$, respectively. 
Then $dim(P_{1}\cap P_{2})\geq n-1$.\\ ii). 
Let $P_{\alpha}$, where $\alpha$ runs in a set $A$, be an arbitrary family of
decomposable $n$-vectors of a linear space
$V$ 
such that every sum
$P_{\alpha_{1}}+P_{\alpha_{2}}$ is also decomposable, and let $D_{\alpha}$ be the
subspaces
spaned by the factors of $P_{\alpha}$, respectively. 
Then either $dim(\cap_{\alpha\in A} D_{\alpha})\geq n-1$ or
$dim(\sum_{\alpha\in A} D_{\alpha})= n+1$. \par
Based on Theorem 2.8, 
Dufour and Zung
prove several linearization theorems, and we refer the reader to \cite{DZ} for these
theorems.
\section{Nambu-Poisson-Lie Groups}
Nambu-Poisson-Lie groups as defined below were discussed in \cite{V3}
and, independently, 
in \cite{GM}, where a complete description of the multiplicative Nambu-Poisson
tensor fields on a Lie group is given.
In this section we reproduce the relevant part of our 
preprint \cite{V3}, and refer the reader
to \cite{GM} for general structural results.

Since Poisson-Lie groups play an important role in Poisson geometry (e.g.,
\cite{V1}), we are motivated to discuss similarly defined Nambu-Poisson-Lie
groups. These cannot be defined by the demand that the multiplication is a
Nambu-Poisson morphism since the direct sum of Nambu-Poisson tensors is
not Nambu-Poisson (it is not decomposable). But, it makes sense to say that
a Nambu-Poisson tensor $P$ endows the Lie group $G$ with the structure of a
{\it Nambu-Poisson-Lie group} if $P$ is a {\it multiplicative tensor field}
i.e. (e.g., \cite{V1}), $\forall g_{1},g_{2}\in G$, one has
$$P_{g_{1}g_{2}}=L_{g_{1}^{*}}P_{g_{2}}+R_{g_{2}^{*}}P_{g_{1}},
\leqno{(3.1)}$$
where $L$ and $R$ denote left and right translations in $G$, respectively.

The multiplicativity of $P$ implies $P_{e}=0$, where $e$ is the unit of $G$.
Moreover, if $G$ is connected, $P$ is multiplicative iff $P_{e}=0$, and the
Lie derivative $L_{X}P$ is a left (right) invariant tensor field whenever
$X$ is left (right) invariant (e.g., \cite{V1}). As an immediate
consequence it follows that
the Nambu-Poisson-Lie group structures on the additive Lie group
${\bf R}^{m}$ are exactly the linear Nambu-Poisson structures
of ${\bf R}^{m}$.

From (3.1), it follows that the set
$$G_{0}:=\{g\in G\;/\;P_{g}=0\}$$
is a closed subgroup. Indeed,
(3.1) shows that if $g_{1}.g_{2}\in
G_{0}$, the product $g_{1}g_{2}\in G_{0}$. Furthermore,
if $g\in G_{0}$, then
$$0=P_{e}=P_{gg^{-1}}=L_{g^{*}}P_{g^{-1}},$$
hence $g^{-1}\in G_{0}$.

In order to give another characterization of Nambu-Poisson-Lie groups, we
generalize a
bracket of $1$-forms, which plays a fundamental role in Poisson geometry (e.g.,
\cite{V1}),
to Nambu-Poisson manifolds.

The natural extension of the bracket of $1$-forms to Nambu-Poisson
structures of order $n$ on $M^{m}$ is defined as follows
$$\{\alpha_{1},\ldots,\alpha_{n}\}=
d(P(\alpha_{(n)}))+\sum_{k=1}^{n}(-1)^{n+k}
i(\sharp_{P}(\alpha_{(n,\hat k)}))d\alpha_{k}
\leqno{(3.2)}$$
$$=\sum_{k=1}^{n}(-1)^{n+k}L_{\sharp_{P}(\alpha_{n,\hat k})}\alpha_{k}
-(n-1)d(P(\alpha_{(n)})),$$
where $\alpha_{k}$ $(k=1,\ldots,n)$ are $1$-forms on $M$.
The equality of the two expressions of
the new bracket follows by using the classical relation
$L_{X}=di(X)+i(X)d$. The bracket (3.2) will be called the
{\it Nambu-Poisson form-bracket}, and we have
\proclaim 3.1 Proposition. The Nambu-Poisson form-bracket satisfies the
following properties: \\i) the form-bracket is totally skew-symmetric;\\
ii) $\forall f_{(n)}\in{\cal F}(M)$, one has
$$\{df_{1},\ldots,df_{n}\}=d\{f_{1},\ldots,f_{n}\};\leqno{(3.3)}$$
iii) for any $1$-forms $\alpha_{(n)}$, and $\forall f\in {\cal F}(M)$ one
has
$$\{f\alpha_{1},\alpha_{2},\ldots,\alpha_{n}\}=
f\{\alpha_{1},\alpha_{2},\ldots,\alpha_{n}\}\leqno{(3.4)}$$
$$+P(df,\alpha_{2},\ldots,\alpha_{n})\alpha_{1}.$$
iv) $\forall f_{(n-1)}\in {\cal F}(M)$ and
for any $1$-form $\alpha$ one has
$$\{df_{1},\ldots,df_{n-1},\alpha\}=L_{X_{f_{(n-1)}}}\alpha.
\leqno{(3.5)}$$
\noindent{\bf Proof.} i) is obvious. ii)
and iii) follow from the first expression
of (3.2).\\
iv) is a consequence of the first expression (3.2) and of
the commutativity of $d$ and $L$.
Q.e.d.

Of course, in view of the skew symmetry formulas corresponding to (3.4),
(3.5) may be used if the factor $f$ and, respectively, the $1$-form
$\alpha$
appear at another factor of the bracket.\vspace{2mm}

It would be nice if the form-bracket would also satisfy the fundamental
identity of Nambu-Poisson brackets. This happens for $n=2$ but, generally,
we only have the following
weaker result
\proclaim 3.2 Proposition. The Hamiltonian vector fields act
as derivations of the
form-bracket by the Lie derivative operation.\par
\noindent{\bf Proof.} Suppose that the required property holds for the
$1$-forms $\alpha_{(n)}$ i.e.,
$$L_{X_{f_{(n-1)}}}\{\alpha_{1},\alpha_{2},\ldots,\alpha_{n}\}=
\sum_{k=1}^{n}\{\alpha_{1},\alpha_{2},\ldots,\alpha_{k-1},
L_{X_{f_{(n-1)}}}\alpha_{k},\ldots,\alpha_{n}\}.\leqno{(3.6)}$$
Then, a straightforward computation which uses (3.4) and (2.5) shows that
$L_{X_{f_{(n-1)}}}$ also acts as a derivation of the bracket
$\{f\alpha_{1},\alpha_{2},\ldots,\alpha_{n}\}$,
$\forall f\in{\cal F}(M)$.

This remark shows that the proposition is true if the
result holds for a bracket of the form
$\{dg_{1},\ldots,dg_{n}\}$, $\forall g_{k}\in{\cal F}(M)$. We see that
this happens by applying (3.3), and the fundamental identity for functions,
since we have
$$L_{X_{f_{(n-1)}}}\{dg_{1},\ldots,dg_{n}\}=
L_{X_{f_{(n-1)}}}d\{g_{1},\ldots,g_{n}\}=
dL_{X_{f_{(n-1)}}}\{g_{1},\ldots,g_{n}\}.$$
Q.e.d.

The relation between (3.6) and the fundamental identity for $1$-forms is
given by (3.5). Moreover, since locally any closed form is an exact form, we
see that the {\it fundamental identity}
$$\{\beta_{1},\ldots,\beta_{n-1},
\{\alpha_{1},\ldots,\alpha_{n}\}\}=
\sum_{k=1}^{n}\{\alpha_{1},\ldots,\alpha_{k-1},\leqno{(3.7)}$$
$$\{\beta_{1},\ldots,\beta_{n-1},\alpha_{k}\},\alpha_{k+1},\ldots,
\alpha_{n}\}$$
{\it holds whenever the} $1$-{\it forms} $\beta$ {\it are closed}.
\vspace{2mm}

Another remark is that, since (3.5) expresses a Lie derivative, it defines a
representation of the Lie algebra ${\cal H}(P)$ of the
real, finite, linear combinations of Hamiltonian vector
fields on the space $\wedge^{1}M$ of the $1$-forms on $M$, and Theorem 3.2
tells us that this representation is by derivations of the form-bracket.
\vspace{2mm}

Now, coming back to Nambu-Poisson-Lie groups, we can extend the following
result of Dazord and Sondaz \cite{DzS}
\proclaim 3.3 Theorem. If $G$ is a connected Lie group endowed with a
Nambu-Poisson tensor field $P$ which vanishes at the unit $e$ of $G$, then
$(G,P)$ is a Nambu-Poisson-Lie group iff the $P$-bracket of any $n$ left
(right) invariant $1$-forms of $G$ is a left (right) invariant $1$-form.
\par
\noindent{\bf Proof.} The same proof as in the Poisson case (e.g.,
\cite{V1}) holds. Namely, the evaluation of the Lie derivative
via (3.2) yields
$$(L_{Y}\{\alpha_{1},\ldots,\alpha_{n}\})(X)=
Y((L_{X}P)(\alpha_{(n)}))\leqno{(3.8)}$$
for any left invariant vector field $X$, right invariant vector field $Y$,
and left invariant $1$-forms $\alpha_{(n)}$. (Same if left and right are
interchanged.) Hence, the condition of the theorem is equivalent with
the fact that $L_{X}P$ is left invariant if $X$ is left invariant. Q.e.d.
\vspace{2mm}

Some other basic properties of Poisson-Lie groups also have a
straightforward generalization.
First of all, since $P_{e}=0$ for a
Nambu-Poisson-Lie group $G$ with unit $e$, and Nambu-Poisson tensor $P$,
the linear approximation of $P$ at $e$ defines a linear Nambu-Poisson structure
on the Lie algebra ${\cal G}$ of $G$, and a dual $n$-Lie algebra structure
on the dual space ${\cal G}^{*}$.
As for $n=2$, a compatibility relation between
the Lie bracket and the linear Nambu-Poisson structure of ${\cal G}$
exists.

First, following \cite{LW},
let us consider the {\it intrinsic derivative} $\pi_{e}:=d_{e}P:
{\cal G}\rightarrow\wedge^{n}{\cal G}$ defined by
$$\pi_{e}(X)(\alpha_{(n)})=(L_{\check X}P)_{e}(\alpha_{(n)})
,\leqno{(3.9)}$$
where $\alpha_{(n)}\in{\cal G}^{*}$, $X\in{\cal G}$, and $\check X$ is any
vector field on {\cal G} with the value $X$ at $e$. Then we have
\proclaim 3.4 Theorem. i). The bracket of the dual $n$-Lie algebra
structure of ${\cal G}^{*}$ is the dual of the mapping $\pi_{e}$, and it
has each of the following expressions
$$[\alpha_{1},\ldots,\alpha_{n}]=
d_{e}(P(\check\alpha_{(n)}))=\pi_{e}^{*}(\alpha_{(n)})
\leqno{(3.10)}$$
$$=\{\bar\alpha_{1},\ldots,\bar\alpha_{n}\}_{e}
=\{\tilde\alpha_{1},\ldots,\tilde\alpha_{n}\}_{e},$$
where $\alpha_{(n)}\in{\cal G}^{*}$, $\check\alpha_{(n)}$ are $1$-forms on $G$
which are equal to $\alpha_{(n)}$ at $e$, and $\bar\alpha_{(n)}$,
$\tilde\alpha_{(n)}$ are the left and right invariant $1$-forms,
respectively, defined by $\alpha_{(n)}$.\\
ii). The mapping $\pi_{e}$ is a $\wedge^{n}{\cal
G}$-valued $1$-cocycle of ${\cal G}$
with respect to the adjoint representation
$$ad\,X(Y_{1}\wedge\ldots\wedge Y_{n})=
\sum_{k=1}^{n}Y_{1}\wedge\ldots Y_{k-1}
\wedge[X,Y_{k}]_{{\cal G}}\wedge Y_{k+1}
\wedge\ldots\wedge Y_{n},$$
($X,Y_{(n)}\in{\cal G}$). \par
\noindent{\bf Proof.} The proofs are exactly the same as for $n=2$; see
\cite{LW} or Chapter 10 of
\cite{V1}. We repeat them briefly here.

i). By the definition of a dual mapping,
and since $P_{e}=0$, we have
$$<\pi_{e}^{*}(\alpha_{(n)}),X>=\pi_{e}(X)(\alpha_{(n)})
=(L_{\check X}P)_{e}(\alpha_{(n)})=X(P(\check
\alpha_{(n)}))$$ $$=<d_{e}(P(\check\alpha_{(n)}),X>,$$
and this differential clearly is the $n$-Lie algebra structure of the
linear approximation of $P$ at $e$. This justifies the first two equality
signs of (3.10). The remaining part of (3.10) follows from:
$$\{\bar\alpha_{1},\ldots,\bar\alpha_{n}\}_{e}(X)\stackrel{(3.2)}{=}
X(P(\bar\alpha_{(n)}))+\sum_{k=1}^{n}
(-1)^{n+k}(d\bar\alpha_{k})_{e}(\sharp_{P}(\alpha_{(n,\hat
k)}),X)$$
$$=X(P(\bar\alpha_{(n)}))-\sum_{k=1}^{n}(-1)^{n+k}
(L_{\tilde X}\bar\alpha_{k})_{e}(\sharp_{P}(\alpha_{(n,\hat
k)}))$$
$$+\sum_{k=1}^{n}\sharp_{P}(\alpha_{(n,\hat
k)})_{e}(\bar\alpha_{k}(\tilde X))=X(P(\bar\alpha_{(n)}),$$
where $\tilde X$ is the right invariant vector field defined by $X$,
and we used the equalities $P_{e}=0$, $L_{\tilde X}\bar\alpha_{k}=0$.
ii). The fact that $\pi_{e}$ is a $1$-cocycle means that we have
$$ad\,X(\pi_{e}(Y))-ad\,Y(\pi_{e}(X))-\pi_{e}([X,Y]_{{\cal G}})=0
,\leqno{(3.11)}$$
where $X,Y\in {\cal G}$.
We always use the notation with bars and tildes for left and right invariant
objects on Lie groups as we did above. Then, it follows that
$$ad\,X(\pi_{e}(Y))=\frac{d}{ds}/_{s=0}Ad\;exp(sX)((L_{\bar Y}P)_{e})=
(L_{\bar X}L_{\bar Y}P)_{e},$$
and (3.11)
is a consequence of this result. Q.e.d.

Now we get the relation announced earlier:
\proclaim 3.5 Corollary. $\forall\alpha_{(n)}\in{\cal G}^{*}$
and $\forall X,Y\in{\cal G}$ the following
relation holds
$$<[\alpha_{1},\ldots,\alpha_{n}],[X,Y]_{{\cal G}}>=
\sum_{k=1}^{n}(<[\alpha_{1},\ldots,\alpha_{k-1},coad_{X}\alpha_{k},
\leqno{(3.12)}$$
$$\alpha_{k+1},\ldots,\alpha_{n}],Y>
-<[\alpha_{1},\ldots,\alpha_{k-1},coad_{Y}\alpha_{k},
\alpha_{k+1},\ldots,\alpha_{n}],X>).$$
\par
\noindent{\bf Proof.} The result is nothing but a reformulation of the
cocycle condition (3.11). Q.e.d. \vspace{2mm}

In agreement with Corollary 3.5,
we will define a {\it Nambu-Poisson-Lie
algebra} as a Lie algebra with a linear Nambu-Poisson structure which
satisfies (3.12).
The question is: given a Nambu-Poisson-Lie algebra
${\cal G}$,
is it possible to integrate it
to a Nambu-Poisson-Lie group? In a forthcomming version of 
\cite{V3} we will show that the general answer is no, even if the 
definition of a Nambu-Poisson-Lie algebra is changed by adding one
more necessary condition which is implied by \cite{GM}.
A corresponding negative example on the unitary Lie algebra $u(2)$
will be quoted later on.

But, some of the results
known for $n=2$ still hold.
If $G$ is connected and simply
connected,
for any $1$-cocycle $\pi_{e}$ as in Theorem 3.4 ii),
there exists a unique multiplicative $n$-vector field $P$ on $G$,
called the {\it integral field} of $\pi_{e}$ such that
$d_{e}P$ is the given cocycle.
Indeed, for the given cocycle $\pi_{e}$, $$\pi_{g}(X_{g}):=Ad\,g(
\pi_{e}(L_{g^{-1*}}X_{g}))\hspace{5mm}(g\in G,\;X_{g}\in T_{g}G)$$
defines a $\wedge^{n}{\cal G}$-valued $1$-form
$\pi$ on $G$ which satisfies the
equivariance condition $L_{g}^{*}\pi=(Ad\,g)\circ\pi$. This implies that
$d\pi=0$, and, since $G$ is connected and simply connected, $\pi=dP$ for a
unique $n$-vector field $P$ on $G$, which can be seen to be
multiplicative \cite{LW} \cite{V1}.
If this field is Nambu-Poisson, we are done. But, this final part
is more complicated than for $n=2$ since it involves the quadratic
identity (2.7), and the non-tensorial differential identity (2.8).
We only have
\proclaim 3.6 Proposition. If ${\cal G}$ is a Nambu-Poisson-Lie algebra of
even order $n$, the integral field $P$
of the dual cocycle $\pi_{e}$ of
the linear Nambu-Poisson structure of ${\cal G}$, on the connected, simply
connected Lie group $G$ which integrates ${\cal G}$,
is a multiplicative generalized Poisson structure on $G$.
\par
\noindent{\bf Proof.} The same proof as for $n=2$ \cite{LW}, \cite{V1}
shows that the Schouten-Nijenhuis bracket $[P,P]=0$. Indeed, since $P$ is
multiplicative, so is $[P,P]$ and, in particular,
$[P,P]_{e}=0$. Furthermore, since $n$ is even, $P_{e}=0$, and using the
coordinate expression of the Schouten-Nijenhuis bracket \cite{V1},
we have
$$d_{e}[P,P](X)=2[P,L_{\check X}P]_{e}=\leqno{(3.13)}$$
$$=\frac{2}{(2n-1)!n!(n-1)!}\delta^{k_{1}...........k_{2n-1}}
_{i_{1}...i_{n}j_{2}...j_{n}}\frac{\partial P^{i_{1}...i_{n}}}
{\partial x^{u}}\left(\xi^{v}\frac{\partial P^{uj_{2}...j_{n}}}
{\partial x^{v}}\right)\frac{\partial}{x^{k_{1}}}\wedge...\wedge\frac{\partial}
{\partial x^{k_{2n-1}}}/_{e},$$
where $X=\xi^{v}(\partial/\partial x^{v})/_{e}$.
Now, $\xi^{v}(\partial P^{j_{1}...j_{n}}/\partial x^{v})$ are the
coordinates of the $n$-vector $(d_{e}P)(X)$ of the linear approximation of
$P$ at $e$. Hence, the result of (3.13) is the
algebraic Schouten-Nijenhuis bracket
$[d_{e}P,d_{e}P]_{{\cal
G}}$ (e.g., \cite{V1}), which is zero by Corollary 2.5. The conclusion is that
$d_{e}[P,P]=0$.
But, a
multiplicative tensor field with a vanishing intrinsic derivative at $e$ is
identically $0$ \cite{LW}, \cite{V1}.
Hence, $[P,P]=0$.
Q.e.d. \vspace{2mm}

Theorem 3.4 also allows us to get a result on subgroups just as in the 
Poisson case. A Lie subgroup $H$ of a 
Nambu-Poisson-Lie group $(G,P)$ will be called 
a {\it Nambu-Poisson-Lie subgroup} 
if $H$ has a (necessarily unique) multiplicative 
Nambu-Poisson tensor $Q$ such that $(H,Q)$ is a 
Nambu-Poisson submanifold of $(G,P)$.
If $H$ is connected, it is a Nambu-Poisson-Lie subgroup of $(G,P)$
iff $Ann({\cal H})$, where ${\cal H}$ is the Lie algebra of $H$,
is an ideal in $({\cal G}^{*},[.,...,.])$. By this we mean that the bracket
(3.10) is in $Ann({\cal H})$ whenever one of the arguments (at least)
is in $Ann({\cal H})$. The proof is the same as for $n=2$ e.g., \cite{V1}.

Furthermore, if $(H,Q)$ is a Nambu-Poisson-Lie subgroup of $(G,P)$,
the homogeneous space $M=G/H$ inherits a Nambu-Poisson structure $S$ of the 
same order as $P,Q$ such that the natural projection
$p:(G,P)\rightarrow(M,S)$ is a Nambu-Poisson morphism. This holds since the 
brackets  $\{f_{1}\circ p,...,f_{n}\circ p\}_{P}$ are
constant along the fibers of $p$, which is easy to check using (3.1).
(E.g., see Proposition 10.30 in \cite{V1} for the case $n=2$.)
Moreover, as a consequence of (3.1), the natural left action
of $G$ on $M$ satisfies the multiplicativity condition
$$S_{g(x)}=\varphi_{g^{*}}(S_{x})+\varphi^{x}_{*}(P_{g}), \leqno{(3.1')}$$
where $\varphi_{g}(x)=\varphi^{x}(g)=g(x)$ for $g\in G,x\in M$, and
$\varphi_{g}:M\rightarrow M,\: \varphi^{*}:G\rightarrow M$.
Accordingly, any action of a 
Nambu-Poisson-Lie group $(G,P)$ on a Nambu-Poisson manifold $(M,S)$ which
satisfies (3.1$'$) will be 
called a {\it Nambu-Poisson action}. If $G$ is connected, one has the
same infinitesimal characteristic properties of Nambu-Poisson actions 
as in the Poisson case e.g., Proposition 10.27 in \cite{V1}.
In particular, that $\forall X\in{\cal G}$,
$L_{X_{M}}S=-[(d_{e}P)(X)]_{M}$, where $e$ is the unit of $G$, and 
the index $M$ denotes the infinitesimal action on $M$.
\vspace{2mm}\\
\indent
Now, we give a number of examples of non commutative Nambu-Poisson-Lie
groups.

A first example is that of the $3$-dimensional solvable Lie group
$$G_{3}:=\{\left( \begin{array}{rrr} x&0&y\\ 0&x&z\\0&0&1 \end{array}
\right)\;\;/x,y,z\in {\bf R},\;x\neq0\}.\leqno{(3.14)}$$
The left invariant forms of this group are $dx/x,dy/x,dz/x$, and if we look
for a Nambu tensor of the form
$$P=f(x)\frac{\partial}{\partial x}\wedge\frac{\partial}{\partial y}\wedge
\frac{\partial}{\partial z}\leqno{(3.15)}$$
such that $\{dx/x,dy/x,dz/x\}$ is left-invariant, and $f(1)=0$, we see that
$f=x(x^{2}-1)/2$ does the job. The corresponding Nambu-%
Poisson-Lie algebra is ${\bf
R}^{3}$ with the linear Nambu structure $x^{1}(\partial/\partial x^{1})
\wedge(\partial/\partial x^{2})\wedge(\partial/\partial x^{3})$.

The next example is that of the {\it generalized Heisenberg group}
$$H(1,p):=\{\left(\begin{array}{rrr}
Id_{p}&X&Z\\ 0&1&y\\0&0&1 \end{array}\right)\},\leqno{(3.16)}$$
where $X=\,^{t}(x_{1}...x_{p})$, $Z=\,^{t}(z_{1}...z_{p})$.
The left invariant $1$-forms of this group are
$$dx_{1},...,dx_{p},dy,dz_{1}
-x_{1}dy,...,dz_{p}-x_{p}dy,\leqno{(3.17)}$$
and $$P=y\frac{\partial}{\partial x_{1}}\wedge
\frac{\partial}{\partial z_{1}}\wedge\frac{\partial}{\partial y}
\leqno{(3.18)}$$
makes $H(1,p)$ into a Nambu-Poisson-Lie group. Indeed, it vanishes at the unit, and
it follows easily that the brackets of the left invariant $1$-forms are
left invariant. The corresponding Nambu-Poisson-Lie
algebra is ${\bf R}^{2p+1}$
with the same Nambu tensor (3.18).

A third example is that of the direct product
$G=H(1,1)\times{\bf R}_{+}$, where ${\bf R}_{+}$ is the multiplicative
group of the positive real numbers $t$. The left invariant $1$-forms
of the group are those given by (3.17), and $dt/t$. The tensor
$$P=t(\ln t)\frac{\partial}{\partial y}\wedge
\frac{\partial}{\partial z}\wedge\frac{\partial}{\partial
t}\leqno{(3.19)}$$
makes $G$ into a Nambu-Poisson-Lie
group for the same reasons as in the previous
examples. The corresponding Nambu-Poisson-Lie
algebra is ${\bf R}^{4}$ with the
linear Nambu structure
$$P=x_{4}\frac{\partial}{\partial x_{2}}\wedge
\frac{\partial}{\partial x_{3}}\wedge\frac{\partial}{\partial
x_{4}}.\leqno{(3.20)}$$

We may notice that if $(G_{1},P)$ is a Nambu-Poisson-Lie
group, and $G_{2}$ is
any other Lie group, $fP$, where $f\in C^{\infty}(G_{2}$, is a
Nambu-Poisson-Lie
structure on $G_{1}\times G_{2}$.

The next example is that of
a  Nambu-Poisson-Lie
algebra.
Consider the unitary Lie algebra $u(2)$ with the basis
$$X_{1}=\frac{\sqrt{-1}}{2}\left(\begin{array}{rr}1&0\\0&1\end{array}\right),\;
X_{2}=\frac{\sqrt{-1}}{2}\left(\begin{array}{rr}0&1\\1&0\end{array}\right),\;$$
$$X_{3}=\frac{1}{2}\left(\begin{array}{rr}0&1\\-1&0\end{array}\right),\;
X_{4}=\frac{\sqrt{-1}}{2}\left(\begin{array}{rr}1&0\\0&-1\end{array}\right).$$
Then, the linear Nambu tensor
$$P=x_{1}\frac{\partial}{\partial x_{2}}\wedge
\frac{\partial}{\partial x_{3}}\wedge\frac{\partial}{\partial
x_{4}}\leqno{(3.21)}$$
yields a Nambu-Poisson-Lie
algebra structure. Indeed, straightforward computations
show that (3.12) is satisfied. In a new version of \cite{V3}, we show
that this structure does not come from a Nambu-Poisson-Lie group. 
Namely, the structure theory of \cite{GM} implies that if (3.21) commes from a
Nambu-Poisson-Lie group structure $\Lambda$ of $U(2)$ then
$$R_{g^{-1}*}\Lambda=\theta(g)\frac{\partial}{\partial x^{2}}
\wedge\frac{\partial}{\partial x^{3}}\wedge
\frac{\partial}{\partial x^{4}}\hspace{5mm}(\forall g\in U(2)),$$
where $\theta$ comes from an additive character of the circle subgroup $S^{1}$ of $U(2)$
hence, $\theta=0$.

In principle, all the Nambu-Poisson-Lie
algebras can be determined from the
Dufour-Zung classification of the linear Nambu structures
\cite{DZ} by looking for
Lie algebras structure constants which, together with the canonical
structures of \cite{DZ}, satisfy the condition (3.12).\vspace{2mm}

Proposition 3.6 might suggest looking for examples of
Nambu-Poisson-Lie groups by first looking for $(2p)$-vector fields
$P$ on a Lie group $G$ which
are multiplicative, and satisfy the Schouten-Nijenhuis bracket condition
$[P,P]=0$. For this purpose, the technique of the {\it
generalized Yang-Baxter
equation}
$$(ad\,X)[{\bf r},{\bf r}]_{{\cal G}}=0\hspace{1cm}(X\in{\cal G},\;{\bf r}\in
\wedge^{2p}{\cal G}),\leqno{(3.22)}$$
used for $n=2$ (e.g., \cite{YKS}, \cite{V1}) may be extended. But, since the
$2p$-vector field to be considered
is $P=\bar{\bf r}-\tilde{\bf r}$
\cite{V1} (remember
that  bar and tilde denote the left and right invariant corresponding
tensor field, respectively), it is not clear whether we can get a
decomposable tensor $P$.\par
On the other hand, we should look for
decomposable solutions of the {\it classical
Yang-Baxter equation}
$$[{\bf r},{\bf r}]=0\hspace{1cm}({\bf r}\in
\wedge^{2p}{\cal G})
\leqno{(3.23)}$$
for another reason too. Namely,
The left (right) invariant field generated by such a solution could give us
left (right) invariant Nambu-Poisson structures on the Lie group $G$.
General questions on left invariant Nambu-Poisson structures 
on Lie groups are studied in
\cite{Nak}.

We end this section by indicating the method used in \cite{GM}
for the construction of the Nambu-Poisson-Lie groups. It consists in looking at 
the sum and the intersection of the subspaces $V_{g}\subseteq{\cal G}$
spaned by the factors of the decomposable $n$-vectors $R_{g^{-1}*}P_{g}$, 
$\forall g\in G$, and showing that these provide an ideal ${\cal H}$
of dimension $n$, $n-1$ or $n+1$ in ${\cal G}$. 
(The result follows from the multiplicativity of $P$ and the use of Lemma 2.9.)
Accordingly,
the multiplicative $n$-vector fields on $G$ are given by acting on the left 
invariant $n$-vector defined by ${\cal H}$ via multiplication by a function $\varphi$,
wedge product by a 
vector field $X$, and interior product by a $1$-form $\alpha$, 
respectively, with well determined properties
described in \cite{GM}. In particular, it turns out that the simple Lie groups
do not admit multiplicative Nambu-Poisson tensors $P$ of order $n\geq3$,
and, if $G=G_{1}\times...\times G_{s}$ is semisimple with simple factors $G_{i}$
$(i=1,...,s)$, the only multiplicative Nambu-Poisson tensors on $G$ are wedge 
products of "contravariant volumes" on a part of the factors with either multiplicative
Poisson bivectors or multiplicative vector fields on other factors.
\section{Questions on quantization}
The quantization of the Nambu-Poisson bracket was considered from the
very first paper \cite{Nb}, and it was discussed by many authors
\cite{{BF},{Tk1},{ChT},{DFS},{DF}}, etc.
This section is not a survey of the quoted references but,
a preliminary discussion
about possible approaches to geometric and
deformation quantization of Nambu-Poisson brackets.

We consider the Kostant-Souriau
geometric quantization \cite{Kt}, \cite{Sr} first.

The prequantization
of a symplectic manifold $M$
is defined by a canonical
lifting of the Hamiltonian vector field $X_{f}$
of an {\it observable} i.e., a function $f\in{\cal F}(M)$, to the total
space of a principal ${\bf C}^{*}$-bundle $p:L^{*}\rightarrow M$
(${\bf C}^{*}={\bf C}\backslash\{0\}$)
or, equivalently, a circle bundle. The lifting is defined by introducing a
Hermitian metric $h$,
and a Hermitian connection $\nabla$ on the associated complex line
bundle $L$. Namely, $\nabla$ decomposes the tangent bundle of $L^{*}$ into
a {\it horizontal} and a {\it vertical part}. The horizontal component of
the prequantization
lift $\hat f$ of $X_{f}$ will be the $\nabla$-horizontal lift of
$X_{f}$, and the vertical component of $\hat f$ will be the infinitesimal
right translation defined on the fibers of $L^{*}$ by the values of $2\pi
\sqrt{-1}f$
along the trajectory of $X_{f}$ which starts at the base point of the
fiber. (The factor $2\pi\sqrt{-1}$ is explained by technical reasons.)
It is shown \cite{Kt} that
$\hat f$ is determined by the conditions
$$p_{*}(\hat f)=X_{f},\;\;\alpha(\hat f)=-2\pi\sqrt{-1}f, \leqno{(4.1)}$$
where $\alpha$ is the connection form of $\nabla$ on $L^{*}$.

Furthermore \cite{Kt},
$\hat f$ can be reinterpreted as a linear operator on the
space $\Gamma(L)$ of the global cross sections of $L$ given by
$$\hat f(\sigma)=\nabla_{X_{f}}\sigma+2\pi\sqrt{-1}f\sigma
\hspace{1cm}(\sigma\in\Gamma(L)),\leqno{(4.2)}$$
and the operator $\hat f$ of (4.2) is called the {\it prequantization} of $f$.
Iff the curvature form $\Omega$ of $\nabla$ satisfies the condition
$$\Omega(X_{f},X_{g})=-2\pi\sqrt{-1}\omega(X_{f},X_{g})
\hspace{1cm}(f,g\in{\cal F}(M)),\leqno{(4.3)}$$
where $\omega$ is the symplectic form of $M$, the prequantization operators
satisfy the celebrated Dirac commutation condition
$$\widehat{\{f,g\}}=[\hat f,\hat g]:=
\hat f\circ\hat g-\hat g\circ\hat f.\leqno{(4.4)}$$

Then, if $L$ is tensorized by the line bundle
$D$ of half-densities (or
half-forms) on $M$, the Hermitian metric $h$ yields a pre-Hilbert
scalar product on
the space $\Gamma_{c}(L\otimes D)$
of the cross sections with compact support of the tensor product $L\otimes D$,
by integration along $M$. The new
{\it prequantization operators}
$\check f:=\hat f\otimes Id+Id\otimes L_{X_{f}}$
are anti-Hermitian with respect to this product
(e.g., \cite{Wd} or the brief survey \cite{V2}).

In the case of a Nambu-Poisson manifold $(M^{m},P)$ of order $n$, a
Hamiltonian vector field is defined by $n-1$ functions $f_{(n-1)}\in
{\cal F}(M)$. Since $\sharp_{P}$ as defined by (2.3) is
multilinear, rather than linear, we may
follow \cite{{TkD},{Nks}}, and
introduce the non associative, non commutative, real
algebra
${\cal O}(M)=\wedge_{{\bf R}}^{n-1}
{\cal F}$
with the product
$$f_{(n-1)}\,\times_{{\cal O}}\,g_{(n-1)}=
\sum_{k=1}^{n-1}g_{1}\wedge\ldots\wedge g_{k-1}\wedge
X_{f_{(n-1)}}g_{k}\wedge g_{k+1}\wedge\ldots\wedge g_{n-1}.\leqno{(4.5)}$$
Then, taking the Hamiltonian vector field extends to a ${\bf R}$-linear mapping
$ham:{\cal O}(M)\rightarrow{\cal H}(P)$
(also denoted by $ham(A)=X_{A}$,
$A\in{\cal O}(M)$), where the Lie algebra
${\cal H}(P)$ is that defined in Section 2. Furthermore,
in view of (2.6), we have
$$ham\,(f_{(n-1)}\,\times_{{\cal O}}\, g_{(n-1)})
=[ham\,(f_{(n-1)}),ham\,(g_{(n-1)})],\leqno{(4.6)}$$
and from (4.6) we get $$X_{(f_{(n-1)}\,
\times_{{\cal O}}\,g_{(n-1)}
+g_{(n-1)}\,\times_{{\cal O}}\,f_{(n-1)})}=0.$$
Accordingly, if we agree to say that $A\in{\cal
O}(M)$ is a {\it Casimir ``function"} of $P$ if $X_{A}
=0$, it follows that
the bracket (4.5) induces a bracket on ${\cal S}(M):={\cal
O}(M)/\{Casimir\:{\rm ``}functions{\rm "}\}$ which makes ${\cal S}(M)$ into a Lie
algebra isomorphic
to ${\cal H}(P)$ \cite{Nks}.

Since $\times_{{\cal O}}$ is not skew symmetric, we consider
the bracket $$[A,B]_{{\cal O}}:=\frac{1}{2}
(A\times_{{\cal O}}B-
B\times_{{\cal O}}A),$$ and
we will say that $({\cal O}(M),[\;\;]_{{\cal O}})$ is the {\it algebra of the
multi-observables} of the Nambu-Poisson manifold $(M,P)$. It may be seen as a
central extension of ${\cal H}(P)$ by the Casimir ``functions" of $P$.
For $n=2$, this is just
the Poisson algebra $({\cal F}(M),\{\})$, and the described construction generalizes the
situation which exists in symplectic and Poisson geometry.

In spite of the fact that $({\cal O}(M),[\;\;]_{{\cal O}})$
is not a Lie algebra for $n>2$, it is
handy to use the terminology of Lie algebra theory whenever the definitions
there naturally extend to our situation. In particular, via the mapping $ham$,
${\cal O}(M)$ has a representation on ${\cal F}(M)$, and we may speak of
${\cal F}(M)$-valued
cochains and their coboundary $\partial$, where $\partial^{2}$ is not
necessarily $0$.

Now, any $1$-cochain $Q:{\cal O}(M)\rightarrow
{\cal F}(M)$ allows
us to define the prequantization of the observable $A\in{\cal O}(M)$ as
being the operator
$$\hat A(\sigma):=\nabla_{X_{A}}\sigma+2\pi\sqrt{-1}Q(A)\sigma,
\leqno{(4.7)}$$
where $\nabla$ and $\sigma$ are as in formula (4.2).

The geometric meaning of $\hat A$ on the principal bundle $L^{*}$
is similar to that of $\hat f$ of (4.1),
and we get
\proclaim 4.1 Theorem. Let $(M,P)$ be a Nambu-Poisson manifold,
and let $Q$ be a $1$-cochain of ${\cal O}(M)$. Then, if
$\partial Q$ satisfies
the condition
$$(\partial Q)(A,B):=X_{A}(Q(B))-X_{B}(Q(A))-
Q([A,B]_{{\cal O}})\leqno{(4.8)}$$
$$=-2\pi\sqrt{-1}\lambda(X_{A},X_{B})
\hspace{1cm}(\forall A,B\in{\cal O}(M)),$$
for some closed $2$-form $\lambda$ which represents an integral
cohomology class of $M$, then there exists a complex line bundle $L$ on
$M$, endowed with a Hermitian metric and connection, such that the
operators (4.7) satisfy the Dirac commutation condition
$$\widehat{[A,B]}_{{\cal O}}=[\hat A,\hat B]:=
\hat A\circ\hat B-\hat B\circ\hat A.\leqno{(4.9)}$$
Conversely, if such a bundle exists, $Q$ satisfies the condition (4.8).
\par
\noindent{\bf Proof.} For an arbitrary $L$ and $\nabla$ as at the beginning
of this section, formula (4.7) leads to the following commutation relation:
$$\widehat{[A,B]}_{{\cal O}}=
\hat A\circ\hat B-\hat B\circ\hat A
+2\pi\sqrt{-1}((\partial Q)(A,B)\leqno{(4.10)}$$
$$+\frac{1}{2\pi\sqrt{-1}}\Omega(X_{A},X_{B})),$$
where $\Omega$ is the curvature $2$-form of $\nabla$. Hence, if (4.9)
holds, we have (4.8) for $\lambda=\Omega$. This is the last assertion of
the theorem. The first part follows from (4.10) again. Indeed, if we have
(4.8) with the integral form $\lambda$, it is well known that there exists
a bundle $L$ with a Hermitian connection $\nabla$ such that $2\pi\sqrt{-1}
\lambda$ is the curvature of $\nabla$ (e.g., \cite{Kt}). Using these $L$
and $\nabla$, we get the desired result. Q.e.d.
\proclaim 4.2 Remark. If $P$ is regular
on an open, dense subset $N$ of $M$, it is enough to quantize the
restriction of the multi-observables to $N$. Thus, we might concentrate on the study
of the quantization of regular Nambu-Poisson manifolds $M$, which have the
simple structure described in Theorem 2.6. (For n=2, this structure is not
so simple, however.) Then, for geometric quantization, it suffices to use
only connections and forms along the leaves of the canonical foliation $S$
of $P$, and replace Theorem 4.1 by the $S$-leafwise version of the
same theorem. \par
For the clarification of this remark,
see the case of the Poisson manifolds in \cite{V1}.
\proclaim 4.3 Remark. The prequantization operators (4.7) act on the
complex, linear space $\Gamma(L)$. But, it is again possible to tensorize
by the halfdensities, and get anti-Hermitian operators on a pre-Hilbert
space $\Gamma(L\otimes D)$
as described earlier for the classical case. \par

A cochain $Q$ which satisfies the hypothesis of
Theorem 4.1, or of its leafwise version, will be called a
{\it quantifier} of the Nambu-Poisson manifold $(M,P)$.
The prequantization
problem reduces to that of
finding good quantifiers but, we have no method to find them.
For $n=2$, the
tautological quantifier $Q(f)=f$ ($f\in {\cal F}$) leads to the classical
geometric quantization.
For $n\geq2$, 
$$Q(A)=\alpha(X_{A})\hspace{1cm}(A\in{\cal O}(M)),\leqno{(4.11)}$$
where $\alpha$ is a $1$-form on $M$, defines a quantifier. For it, we have
$(\partial Q)(A,B)=
d\alpha(X_{A},X_{B})$, and we may use 
the trivial bundle $L$ with the connection defined by the connection form
$-2\pi\sqrt{-1}\alpha$ as a prequantization bundle. This yields
$\hat A=X_{A}$, which is a trivial quantization, while
what we need is a non trivial quantization.

It is to be noted that if $Q$  is a $1$-cocycle i.e., $\partial Q=0$, we obtain a
prequantization which satisfies the Dirac condition on the trivial complex line bundle
over $M$.

Following is an exaple of a $1$-cochain on ${\cal O}(M)$ which shows
the basic difficulty in finding a quantifier.
Namely,
let $Y_{1},
\ldots,Y_{n-2}$ be arbitrary vector fields on $M$, and put
$$Q(f_{(n-1)})
=det(f_{(n-1)},Y_{1}f_{(n-1)},\ldots,Y_{n-2}f_{(n-1)}),\leqno{(4.12)}$$
where the $(n-1)$-dimensional vectors included are the columns of the
determinant.
Then, the
properties of a determinant show that $Q$ extends to a well defined
$1$-cochain of ${\cal O}$, and we get
$$(\partial Q)(f_{(n-1)},g_{(n-1)})=
\sum_{k=1}^{n-1}(-1)^{k}[det(U_{g};f_{(n-1,\hat k)},Y_{1}f_{(n-1,\hat k)},\ldots,
Y_{n-2}f_{(n-1,\hat k)})$$
$$-det(U_{f};g_{(n-1,\hat k)},Y_{1}g_{(n-1,\hat k)},\ldots,
Y_{n-2}g_{(n-1,\hat k)})]
+Q([f_{(n-1)},g_{(n-1)}]_{\cal O}),$$ where,
$U_{f}$ is the operation of adding
at the top of each column of the remaining
matrix
$0$ on the first column, and
$[X_{f_{(n-1)}},Y_{k-1}]g_{k}$
on the $k^{{\rm th}}$ column, and $U_{g}$ is similar but with the roles of
$f$ and $g$ interchanged.
The $1$-cochain $Q$ of (4.12) is a generalization of the tautological
quantifier of the Poisson case but, it is not a quantifier for $n\geq3$
since $(\partial Q)(f_{(n-1)},g_{(n-1)})$ depends on the functions and not
just on the corresponding Hamiltonian vector fields.

One possible way to avoid this difficulty is restrict prequantization to a
subalgeba of ${\cal O}(M)$, in the spirit of the second step, quantization,
in classical geometric quantization theory.
For instance, let ${\cal S}$ be the subalgebra of the
elements $A\in{\cal O}(M)$ such that $[Y_{i},X_{A}]=0$, $\forall
i=1,...,n-2$, and let ${\cal C}$ be an Abelian subalgebra of ${\cal S}$.
Then, the expression of $\partial Q$ given above shows that the restriction
of $Q$ to ${\cal C}$ is a cocycle on this latter subalgebra,
and $Q/_{{\cal C}}$ allows us to do
geometric prequantization on the trivial complex line bundle. 
A second way out of the mentioned
difficulty would be to conveniently change the definition of the bracket
$[\;]_{{\cal O}}$.\vspace{2mm}

Now, let us refer to deformation quantization.
It was shown by Dito, Flato, Sternheimer and
Takhtajan \cite{{DF},{DFS}} that the deformation quantization of Nambu-Poisson
brackets of order $n\geq3$ should be done via a preliminary
{\it Abelian deformation} of the usual product
of functions. Again, it suffices to study only
regular Nambu-Poisson brackets i.e., brackets defined by a Jacobian
determinant (see Section 2). The basic remark \cite{DFS} is that a Jacobian
determinant defines a Nambu-Poisson bracket because the usual product of
functions satisfies the following properties: a) associativity, b)
commutativity, c) distributivity, d) the Leibniz rule of derivation. Hence,
any deformation of the usual product which continues to satisfy a), b), c),
d) allows us to define a deformed Jacobian which is a
Nambu-Poisson bracket on the deformed algebra of $C^{\infty}$-functions on
$M$. (In \cite{SaV} the authors claim the non-existence of a Nambu-Poisson
deformation quantization on $C^{\infty}(M)$ itself.)

In a different formulation, let $(M,P)$ be a regular Nambu-Poisson manifold
of order $n\geq3$, which has the bracket defined by formula (2.22). Assume
that there exists an embedding of complex, linear spaces
$$\iota:{\cal F}(M,{\bf
C})
\rightarrow {\cal A}_{\nu}
:={\cal F}(M,{\bf C})[[\nu]],\leqno{(4.13)}$$
where ${\cal F}(M,{\bf C})$ is the algebra of
complex valued, differentiable functions on $M$,
$\nu$ is a parameter, and
${\cal A}_{\nu}$ is the linear space
of formal power series, endowed with a product $*_{\nu}$ which makes it
an associative,
commutative algebra. The product $*_{\nu}$ of ${\cal A}_{\nu}$ is called an
{\it Abelian product deformation}, and $\forall f,g\in{\cal F}(M,{\bf C})$ one
defines the {\it star product} $f*_{\nu}g:=\iota(f)*_{\nu}\iota(g)$.
Assume also
that the Lie algebra $\chi(M)$ of the vector fields on $M$ has a
representation $\rho$ by derivations
of $({\cal A}_{\nu},*_{\nu})$
Then, we can define ${\cal A}_{\nu}$-valued forms, and their
$*_{\nu}$-exterior product and $\rho$-exterior differential $d_{\rho}$ by
extending the classical definitions. Now, if we write (2.22) for these new
operations,  we get
$$d_{\rho S}(\iota f_{1})\wedge\ldots\wedge d_{\rho S}(\iota f_{n})=
\{f_{1},\ldots,f_{n}\}_{\nu}(\iota\omega),\leqno{(4.14)}$$
where $\iota\omega$ is defined by
$$(\iota\omega)(X_{1},...,X_{n})=\iota(\omega(X_{1},...,X_{n})).$$
The bracket $\{f_{1},\ldots,f_{n}\}_{\nu}$ is the {\it quantum
deformation} of $\{f_{1},\ldots,f_{n}\}$, and it satisfies the properties
of a Nambu-Poisson bracket (i.e., i), ii), iii) of Section 2).

In \cite{{DF},{DFS}}, the
authors propose a construction of an Abelian product
$*_{\nu}$, which leads to a quantum deformation of a Nambu-Poisson bracket
called
{\it Zariski quantization}. For this theory we refer the reader to the
quoted original papers.

Here, we
modify a construction used in symplectic deformation
quantization
\cite{{Lich},{Fds}}
in order to get a
deformation of the Nambu-Poisson bracket if the algebra
${\cal F}(M,{\bf C})$ is embedded into a larger algebra
$\tilde{\cal F}(M,{\bf C})$ first, and the space $\tilde{\cal A}_{\nu}
:=\tilde{\cal F}(M,{\bf C})[[\nu]]$ of
formal power series is used. This construction is not an answer to
the deformation quantization problem since the
obtained 
$*_{\nu}$-product of functions is a power series with coefficients
which may not be functions. It is an example of a general,
commutative, product deformation process,
associated with a fixed Riemannian metric $g$ on the
Nambu-Poisson manifold $(M,P)$.

Let us introduce the associative, commutative
algebra
$$\tilde{\cal F}(M,{\bf C})=\oplus_{i=0}^{\infty}\Gamma\odot^{i}T_{c}^{*}M
,\leqno{(4.15)}$$
where $\Gamma\odot^{i}T_{c}^{*}M$ is the space of symmetric, $i$-covariant,
complex tensor fields, any
particular element of $\tilde{\cal F}(M,{\bf C})$ consists of a finite sum
of terms,
and the product
in the algebra (4.15)
is the symmetric tensor product $\odot$.
As in \cite{Lich}, \cite{Fds}, we define a {\it Weyl-Moyal
product} of power series
$$a_{u}=\sum_{k=0}^{\infty}\sum_{i}\nu^{k}a_{(u)k}^{i}
\in\tilde{\cal A}_{\nu}
=\tilde{\cal F}(M,{\bf C})[[\nu]]\hspace{5mm}(u=1,2)\leqno{(4.16)}$$
by the formula
$$a_{1}*_{\nu}a_{2}=\sum_{p=0}^{\infty}\frac{\nu^{p}}{p!}
(\partial^{p}a_{1},\partial^{p}a_{2})_{g},\leqno{(4.17)}$$
where the {\it algebraic derivative} $\partial$ is defined on each term of
the series (4.16) as the operator $\partial:\odot^{i} T_{c}^{*}M\rightarrow
Hom(\odot^{i-1} T_{c}^{*}M,T_{c}^{*}M)$ given by
$$(\partial t)(X_{1},\ldots,X_{i-1})(Y):=t(Y,X_{1},\ldots,X_{i-1}),$$
$t\in \odot^{i} T_{c}^{*}M$, and
all the arguments are tangent vectors. Of course, $\partial^{p}$ is the
iteration of $\partial$. Finally, $(\;,\;)_{g}$
is the scalar product induced by $g$.
(In the symplectic case,
there was a symplectic scalar product instead.) Here, the
symmetry of $g$ ensures that formula (4.17) defines the structure
of an associative,
commutative algebra on $\tilde{\cal A}_{\nu}$.

Furthermore,
the action of any vector field $X$ on $M$ as a directional
derivative of functions extends to $\tilde{\cal A}_{\nu}$ by means of
the covariant derivative $\nabla_{X}$ of the tensor fields
with respect to the Riemannian connection of $g$. This action is a
representation $\rho$ by derivations. Accordingly, the
Nambu-Poisson bracket $P$ gets deformed to a Nambu-Poisson bracket on
$\tilde{\cal A}_{\nu}$.

Now, we have to consider an embedding $\iota:{\cal F}(M,{\bf
C})
\rightarrow\tilde {\cal A}_{\nu}$ e.g., the
{\it gradiental
deformation} $$\iota(f)=f+\sum_{i=1}^{\infty}\frac{\nu^{i}}{i!}(\odot^{i}df)
\hspace{1cm}(f\in{\cal F}(M,{\bf C})),
\leqno{(4.18)}$$
then put
$$f*_{\nu}k:=(\iota f)*_{\nu}(\iota k)\hspace{1cm}(f,k\in{\cal F}(M,{\bf C})),$$
as given by (4,17).

Then, we might look at the
``semi-classical
approximation" i.e., take only the term $i=1$ in (4.18). This yields
$$f*_{\nu}k=fk+\nu(fdk+kdf)+\nu^{2}(df\odot dk)
+\nu^{3}(df,dk)_{g}.\leqno{(4.19)}$$
The result
is a polynomial deformation of the
product which has symmetric tensor fields as
coefficients. This product is commutative, and associative, since it is a
restriction of $*_{\nu}$. Then, if we define
$\{f_{1},\ldots,f_{n}\}_{\nu}$
by formula (4.14) interpreted on
$\tilde {\cal A}_{\nu}$, we get a polynomial
deformation of the $P$-bracket of functions $\{f_{1},\ldots,f_{n}\}$, with
symmetric tensor fields as coefficients, which satisfies all the
axioms of a Nambu-Poisson bracket.

If, instead of (4.19), we take the star product
$$f*_{\lambda}k:=fk+\lambda(df,dk)_{g},\leqno{(4.20)}$$
where $\lambda=\nu^{3}$ is the new deformation parameter,
this product is also commutative, but it is associative only in the
semi-classical approximation i.e., up to terms in $\lambda^{k}$ with
$k\geq2$. Indeed, (4.20) implies
$$(f*_{\lambda}k)*_{\lambda}l=fkl +\lambda[f(dk,dl)_{g}
+k(dl,df)_{g}\leqno{(4.21)}$$
$$+l(df,dk)_{g}]
+\lambda^{2}(dl,d(df,dk)_{g})_{g},$$
which justifies the previous assertion.
 \vspace*{1cm}
{\small Department of Mathematics, \\}
{\small University of Haifa,\\}
{\small E-mail: vaisman@math.haifa.ac.il\\}
\end{document}